\documentclass[11pt]{article}
\usepackage{latexsym,amssymb,amsmath,amsfonts,amsthm}

\newtheorem{thm}{\textbf Theorem}[section]
\newtheorem{lem}{\textbf Lemma}[section]
\newtheorem{rem}{\textbf Remark}[section]
\newtheorem{cor}{\textbf Corollary}[section]

\newtheorem{prop}{\textbf Proposition}[section]
\newtheorem{defin}{\textbf Definition}[section]

\newcommand{\md}{\mbox{d}}
\newcommand{\be}{\begin{eqnarray}}
\newcommand{\ee}{\end{eqnarray}}
\newcommand{\mr}{\mathbb{R}}
\newcommand{\mc}{\mathcal}
\newcommand{\mb}{\mathbb}
\newcommand{\mx}{\mbox}
\newcommand{\ep}{\varepsilon}
\newcommand{\tr}{\triangle}
\newcommand{\pt}{\partial}

\begin{document}
\begin{titlepage}
\title{\bf On well-posedness of the Cauchy problem for MHD system in Besov spaces}
\author{Changxing Miao$^1$, Baoquan Yuan$^2$ \\
       \\
       $^1$ Institute of Applied Physics and Computational Mathematics,\\
        P.O. Box 8009, Beijing 100088, P.R. China.\\
        (miao\_{}changxing@iapcm.ac.cn)\\
        \\ $^2$ College of Mathematics and Informatics, Henan Polytechnic University,
        \\  Jiaozuo City, Henan Province, 454000, P.R. China.\\
          (bqyuan@hpu.edu.cn)}
\date{}
\end{titlepage}
\maketitle
\begin{abstract}This paper is devoted to the study of the Cauchy problem of incompressible
magneto-hydrodynamics system in framework of Besov spaces. In the
case of  spatial dimension $n\ge 3$ we establish the global
well-posedness of  the Cauchy problem of incompressible
magneto-hydrodynamics system  for small data and the local one for
large data in Besov space $\dot{B}^{\frac np-1}_{p,r}(\mr^n)$,
$1\le p<\infty$ and $1\le r\le\infty$. Meanwhile, we also prove
the weak-strong uniqueness of solutions with data in
$\dot{B}^{\frac np-1}_{p,r}(\mr^n)\cap L^2(\mr^n)$ for $\frac
n{2p}+\frac2r>1$. In case of $n=2$, we establish the global
well-posedness of solutions  for large initial data in homogeneous
Besov space $\dot{B}^{\frac2p-1}_{p,r}(\mr^2)$ for $2< p<\infty$
and $1\le r<\infty$.
 \vskip0.1in


\end{abstract}

\vspace{.2in} {\bf Key words:}\quad Incompressible
magneto-hydrodynamics system, homogeneous Besov space,
well-posedness, weak-strong uniqueness.




\section{Introduction}
\setcounter{equation}{0} \setcounter{equation}{0} In this paper we
consider the $n$-dimensional incompressible magneto-hydrodynamics
(MHD) system
\begin{eqnarray}\label{1.1}
u_t-\triangle u+(u\cdot \nabla)u-(b\cdot\nabla )b-\nabla p=0\\
\label{1.2} b_t-\triangle b+(u\cdot\nabla)b-(b\cdot\nabla)u=0\\
\label{1.3} \mbox{div}u=0,\ \ \mbox{div}b=0
\end{eqnarray}
with initial data
\begin{eqnarray}\label{1.4}
u(0,x)=u_0(x),\\
\label{1.5} b(0,x)=b_0(x).
\end{eqnarray}
where $x\in \mathbb{R}^n$, $t>0$. Here $u=u(t,x)=(u_1(t,x),\cdots,
u_n(t,x))$, $b=b(t,x)=(b_1(t,x),\cdots, b_n(t,x))$ and $p=p(t,x)$
are non-dimensional quantities corresponding to the flow velocity,
the magnetic field and the pressure at the point $(t,x)$, and
$u_0(x)$ and $b_0(x)$ are the initial velocity and initial
magnetic field satisfying div$u_0$=0, div$b_0$=0, respectively.
For simplicity, we have  included the quantity $\frac12|b(t,x)|^2$
into $p(t,x)$ and we set the Reynolds number, the magnetic
Reynolds number, and the corresponding coefficients to be equal to
$1$.

It is well known that for any initial data $(u_0,b_0)\in
L^2(\mr^n)$ with $n\ge 2$, the MHD equations
(\ref{1.1})-(\ref{1.5}) have been shown to possess at least one
global $L^2$ weak solution $(u(t,x),b(t,x))\in
C_b([0,T];L^2(\mr^2))\cap
L^2((0,T]);\dot{H}^1(\mr^2))$ for any $T>0$ such that
\be\label{1.6}
\|(u,b)\|^2_{L^2(\mr^2)}+2\int^t_0\|(\nabla u(s),\nabla
b(s))\|^2_{L^2(\mr^2)}\md s\le\|(u_0,b_0)\|^2_{L^2(\mr^2)},
 \ee
  but the uniqueness and regularity
remain open besides the case of $n=2$, \cite{D-L,S-T}.
Usually, we  define a Leray weak solution by
 any $L^2$ weak solution $(u,b)$ to the MHD
(\ref{1.1})-(\ref{1.5}), i.e. which  satisfies the MHD equations in
distribution sense, and satisfying the energy estimate (\ref{1.6}).

When
$n=2$, for initial data $(u_0(x),b_0(x))\in L^2(\mr^2)$ there
exists a unique global solution to MHD system
(\ref{1.1})-(\ref{1.3})  with   $(u(t,x),b(t,x))\in
C_b([0,\infty);L^2(\mr^2))\cap
L^2((0,\infty);\dot{H}^1(\mr^2))\cap
C^{\infty}((0,\infty)\times\mathbb{R}^2)$, where $C_b(I)$ denotes
the space of bounded and continuous functions on $I$
\cite{D-L,S-T}. Note that the coupled relation between equations
(\ref{1.1}) and (\ref{1.2}) as well as the relation
$$((b\cdot\nabla)b,u)+((b\cdot\nabla)u,b)=0,\ \mx{ for any }0\le
t<\infty, $$
where $(\cdot,\cdot)$ stands for the inner product in  $L^2$ with respect to
the spatial variables. It follows that the solution $(u,b)$
satisfies the energy equality:
 \be \label{1.7}
\|(u,b)\|^2_{L^2(\mr^2)}+2\int^t_0\|(\nabla u(s),\nabla
b(s))\|^2_{L^2(\mr^2)}\md s=\|(u_0,b_0)\|^2_{L^2(\mr^2)},
 \ee
for any $0\le t<\infty$.

The purpose of this paper can be divided into two aspects. At first,
we prove that for initial data $(u_0,b_0)\in
\dot{B}^{n/p-1}_{p,r}(\mr^n)$, $1\le r\le \infty$, $1\le
p<\infty$, the Cauchy problem (\ref{1.1})-(\ref{1.5}) has the unique local
strong solution
 or global strong small solution in Besov space $\dot{B}^{n/p-1}_{p,r}(\mr^n)$.
If we further assume that the data
$(u_0,b_0)$ is in $L^2(\mr^n)$, the above solution coincides with any
Leray weak solution associated with $(u_0,b_0)$. In fact, we shall establish the stability
result of the Leray weak
solution and strong solution  in Section 3 which implies  the weak and strong
uniqueness.

\begin{thm}\label{thm1.1}
Let $(u_0,b_0)\in \dot{B}^{n/p-1}_{p,r}(\mr^n)$, $1\le p<\infty$,
$1\le r\le\infty$ , $2< q\le \infty$ and $\mx{\rm div}u_0=\mx{\rm div}b_0=0$.

(i) For $1\le r\le \infty$, there exists $\ep_0>0$ such that if
$\|(u_0,b_0)\|_{\dot{B}^{n/p-1}_{p,r}}< \ep_0$, then
(\ref{1.1})-(\ref{1.5}) has a unique solution $(u,b)$ satisfying
 \be\label{1.8}
(u,b)\in C_b(\mr^+;\dot{B}^{n/p-1}_{p,r})\cap
\widetilde{L}^q(\mr^+;\dot{B}^{s_p+2/q}_{p,r}(\mr^n)), \quad  r<\infty,
 \ee
 or
\be\label{1.9}
(u,b)\in C_*(\mr^+;\dot{B}^{n/p-1}_{p,\infty})\cap
\widetilde{L}^q(\mr^+;\dot{B}^{s_p+2/q}_{p,\infty}(\mr^n)), \quad  r=\infty,
 \ee
where $s_p=\frac np-1>1-\frac4q$ is a real number.

(ii) For $1\le r< \infty$, there exists a time $T$ and a unique
local solution  $(u(t,x),b(t,x))$ to the system (\ref{1.1})-(\ref{1.5})
 such that
 \be\label{1.10}
(u,b)\in C_b([0,T];\dot{B}^{n/p-1}_{p,r})\cap
\widetilde{L}^q([0,T];\dot{B}^{\frac np+\frac2q-1}_{p,r}(\mr^n)),\quad r<\infty,
 \ee
or
\be\label{1.11}
(u,b)\in C_*([0,T];\dot{B}^{n/p-1}_{p,\infty})\cap
\widetilde{L}^q([0,T];\dot{B}^{s_p+2/q}_{p,\infty}(\mr^n)), \quad  r=\infty,
 \ee
where $p,\ q$ satisfying $\frac n{2p}+\frac 2q>1$,  $C_*$ denote
the continuity in $t=0$ with respect to time $t$  in weak star sense,
$\widetilde{L}^q([0,T];\dot{B}^{\frac np+\frac2q-1}_{p,r}(\mr^n))$
denotes  the mixed space-time space defined by Littlewood-Paley theory,
please refer to Section 2 for  details.
\end{thm}

\begin{thm} \label{thm1.2}
Let $(u_0,b_0)\in\dot{B}^{n/p-1}_{p,r}(\mr^n)\cap L^2(\mr^n)$ be a
divergence free datum. Assume $1\le p<\infty$ and $2<r<\infty$
such that $\frac n{2p}+\frac2r>1$. Let $(u,b)\in
C([0,T];\dot{B}^{n/p-1}_{p,r}(\mr^n))\cap
L^\infty(\mr^+;L^2(\mr^n))\cap L^2(\mr^+;\dot{H}^1(\mr^n))$ be the
unique solution associated with $(u_0,b_0)$. Then all Leray
solutions associated with $(u_0,b_0)$ coincide with $(u,b)$ on the
interval $[0,T]$.
\end{thm}

Secondly, we shall establish the global well-posedness for the Cauchy problem of the MHD
system (\ref{1.1})-(\ref{1.5}) for data in larger space than
$L^2(\mr^2)$ space, i.e. the homogeneous Besov space
$\dot{B}^{2/p-1}_{p,r}(\mr^2)$ for $2<p<\infty$ and $1\le r<\infty$.
 Let us give some rough analysis. If $1\le
p<2$ and $1\le r\le \infty$ or $p=2$ and $1\le r\le 2$, the global well-posedness is trivial
because of the embedding relation
$\dot{B}^{2/p-1}_{p,r}(\mr^2)\hookrightarrow L^2(\mr^2)$;  The case
$2\le p<\infty$ and $1\le r\le 2$  can be deduced into the case
$2\le p<\infty$ and $2<r<\infty$ because of Sobolev embedding
$\dot{B}^{2/p-1}_{p,r_1}(\mr^2)\hookrightarrow \dot{B}^{2/p-1}_{p,r_2}(\mr^2)$ with $r_1\le r_2$.
 An interesting
question  is whether the MHD system (\ref{1.1})-(\ref{1.5}) is
global well-posedness for arbitrary data in the Besov space
$\dot{B}^{2/p-1}_{p,r}(\mr^2)$ for $2\le p<\infty$, $r=\infty$.

\begin{thm}\label{thm1.3}
Let $(u_0(x),b_0(x))\in \dot{B}^{2/p-1}_{p,r}(\mr^2)$ be
divergence free vector field. Assume that $2\le p<\infty$ and $1\le
r<\infty$. Then there exists a unique solution to the MHD system
(\ref{1.1})-(\ref{1.5}) such that $(u,b)\in
C([0,\infty);\dot{B}^{2/r-1}_{p,r}(\mr^2))$. Moreover, if $p,\ r$
satisfy also $\frac2p+\frac2r>1$ and $1\le r<\infty$, the
following estimate holds:
 \be\label{1.12}
\|(u,b)\|_{\dot{B}^{2/p-1}_{p,r}}\le
C\|(u_0,b_0)\|^{1+\beta}_{\dot{B}^{2/p-1}_{p,r}}
 \ee
for any $t\ge 0$, where $\beta>\frac p2$.
\end{thm}

  From the above discussion, it is sufficient to  prove  the case $2\le p<\infty$ and $2<
r<\infty$ in Theorem \ref{thm1.3}. Since  $(u_0,b_0)\in
C([0,\infty);\dot{B}^{2/r-1}_{p,r}(\mr^2))$ has infinite energy, so we have to use Caldr\'{o}n's argument
\cite{Calderon, G-P} and  perform a interpolation  between the $L^2$-strong solution and
the solution in $C([0,\infty);\dot{B}^{2/{\bar r}-1}_{\bar p,\bar r}(\mr^2))$ with
$p<\bar{p}<\infty$ and $r<\bar{r}<\infty$. In detail, let us decompose data
 \be\label{1.13}
 (u_0(x),b_0(x))=(v_0(x),g_0(x))+(w_0(x),h_0(x)),
 \ee
with $(v_0,g_0)\in L^2(\mr^n)$ and $(w_0,h_0)\in
\dot{B}^{2/\bar{p}-1}_{\bar{p},\bar{r}}(\mr^2)$ for some
$p<\bar{p}<\infty$ and $r<\bar{r}<\infty$ with small norm. The
corresponding solutions are denoted by $(v(t,x),g(t,x))$ and
$(w(t,x),h(t,x))$, where the solutions $(w,h)$ satisfies the MHD
system and $(v,g)$ satisfies MHD-like equations. The global
existence of solution $(w,h)$ in the Besov space
$L^q((0,\infty);\dot{B}^{n/p+2/q-1}_{p,r}(\mr^n))$, $1\le
p<\infty$, $2<q\le\infty$ and $1\le r\le\infty$ for $n\ge 2$ can
be generally proved. The MHD-like system is locally solved, then by
the energy inequality we prove $(v,g)$ is
global solvable for $n=2$. The idea  comes from I. Gallagher and F. Planchon \cite{G-P} who deal
with the Navier-Stokes equations, however  we have give a
different proof for the strong solutions to the MHD system
(\ref{1.1})-(\ref{1.5}) on the mixed time-space Besov spaces
$\widetilde{L}^q([0,T];\dot{B}^{\frac
np+\frac2q-1}_{p,r}(\mr^n))$.

The remaining parts of the present paper are organized as follows.
Section 2 gives some definitions and preliminary tools. In Section 3 we establish some linear estimates
and bilinear estimates of the solution
in framework of mixed space-time Besov space by Fourier localization
and Bony's para-product decomposition,  and by which we complete the proof
 of Theorem\ref{thm1.1} and Theorem\ref{thm1.2}.
Theorem\ref{thm1.3} will be proved in Section 4 by Caldr\'{o}n's argument in conjunction with
the real interpolation method.

We conclude this section by introducing some notations.
Denote by ${\cal S}(\mathbb{R}^n)$ and ${\cal S}'(\mathbb{R}^n)$
the Schwartz space and the Schwartz distribution space, respectively. For any
interval $I \subset \mathbb{R}$ and any Banach space X we denote by C(I ;X) the space of strongly continuous
functions from $I$ to $X$,  and  by $\mc{C}_\sigma(I;B)$  the time-weighted
space-time Banach space  as follows
$$\mc{C}_\sigma(I;X)=\Big\{f\in C(I;B): \; \|f;\mc{C}_\sigma(I;X)\|
  =\sup_{t\in I}t^{\frac1\sigma}\|f\|_X<\infty\Big\}.$$
we denote by $L^q (I ;X)$  and $L^{q_1,q_2}(I; X)$ the space of strongly measurable functions
from $I$ to $X$ with $\|u(\cdot);X\| \in L^q (I )$  and $\|u(\cdot);X\| \in L^{q_1,q_2 }$, respectively.
 $ L^{q_1,q_2}$ denotes usual  Lorentz space,
please refer to \cite{B-L, M1, S-W} for details .

\vskip0.2cm

\noindent{\bf Notation:} Throughout the paper, $C$ stands for   a generic  constant.
We  will use the notation $A\lesssim B$ to denote the relation  $A\le CB$ and
the notation $A\approx B$ to denote the relations  $A\lesssim B$ and $B\lesssim A$.
Further, $\|\cdot\|_{p}$ denotes the norm of the Lebesgue space $L^p$ and
$\|(f_1, f_2, \cdots, f_n)\|_{X}^a$ denotes $\|f_1\|_{X}^a+  \cdots+\|f_n\|_{X}^a$.
The time interval $I$ may be either $[0,T)$
for any $T>0$ or $[0,\infty)$.


\section{Preliminary}
\setcounter{equation}{0}

In this section we first introduce Littlewood-Paley decomposition
and the definition of Besov spaces. Given $f(x)\in \mc{S}(\mr^n)$,
 define the
Fourier transform as
 \be
\hat{f}(\xi)=\mc{F}f(\xi)=(2\pi)^{-n/2}\int_{\mr^n}\mx{e}^{-ix\cdot\xi}f(x)\md
x,
 \ee
and its inverse Fourier transform:
 \be
\check{f}(x)=\mc{F}^{-1}f(x)=(2\pi)^{-n/2}\int_{\mr^n}\mx{e}^{ix\cdot\xi}f(\xi)\md\xi.
 \ee
 Choose two nonnegative radial functions $\chi$, $\varphi \in {\cal
S}(\mathbb{R}^n)$ supported respectively in ${\cal B}=\{\xi\in\mathbb{R}^n,\,
|\xi|\le\frac{4}{3}\}$ and ${\cal C}=\{\xi\in\mathbb{R}^n,\,
\frac{3}{4}\le|\xi|\le\frac{8}{3}\}$ such that
\be
\chi(\xi)+\sum_{j\ge0}\varphi(2^{-j}\xi)=1,\quad\xi\in\mathbb{R}^n,\\
\sum_{j\in\mathbb{Z}}\varphi(2^{-j}\xi)=1,\quad\xi\in\mathbb{R}^n\backslash \{0\}.
\ee
Set $\varphi_j(\xi)=\varphi(2^{-j}\xi)$ and
let $h={\cal F}^{-1}\varphi$ and $\tilde{h}={\cal F}^{-1}\chi$.
Define the frequency localization operators:
\be
&&\Delta_jf=\varphi(2^{-j}D)f=2^{nj}\int_{\mathbb{R}^n}h(2^jy)f(x-y)dy, \\
&&S_jf=\sum_{k\le j-1}\Delta_kf=\chi(2^{-j}D)f=2^{nj}\int_{\mathbb{R}^n}\tilde{h}(2^jy)f(x-y)dy.
\ee
Formally, $\Delta_j=S_{
j}-S_{j-1}$  is a frequency projection into the annulus
$\{|\xi|\approx 2^j\}$, and  $S_j$ is a frequency projection into the
ball $\{|\xi|\lesssim 2^j\}$. One easily verifies that
with the above choice of $\varphi$
\begin{eqnarray}\label{2.1}
\Delta_j\Delta_kf\equiv0\quad i\!f\quad|j-k|\ge 2\quad and
\quad \Delta_j(S_{k-1}f\Delta_k
f)\equiv0\quad i\!f\quad|j-k|\ge 5.
\end{eqnarray}
We now introduce the following definition of Besov spaces.
\begin{defin}\label{Def2.1}Let
$s\in \mathbb{R}, 1\le p,q\le\infty$. The homogenous Besov space $\dot
{B}^s_{p,q}$ is defined by
$$\dot {B}^s_{p,q}=\{f\in {\cal Z}'(\mathbb{R}^n):   \|f\|_{\dot
{B}^s_{p,q}}<\infty\},$$
\end{defin}
\noindent where
$$\|f\|_{\dot{B}^s_{p,q}}=\left\{\begin{array}{l}
\displaystyle\bigg(\sum_{j\in\mathbb{Z}}2^{jsq}\|\Delta_j f\|_p^q\bigg)^{\frac 1
q},\quad \hbox{for}\quad q<\infty,\\
\displaystyle\sup_{j\in \mathbb{Z}}2^{js}\|\Delta_jf\|_p, \quad \hbox{ for}
\quad q=\infty,
\end{array}\right.
$$
and ${\cal Z}'(\mathbb{R}^n)$ can be identified by the quotient
space ${\cal S}'/{\cal P}$ with the  space ${\cal
P}$ of polynomials.
\begin{defin}\label{Def2.2}Let
$s\in \mathbb{R}, 1\le p,q\le\infty$. The inhomogeneous Besov space $
{B}^s_{p,q}$ is defined by
$${B}^s_{p,q}=\{f\in {\cal S}'(\mathbb{R}^n): \|f\|_{
{B}^s_{p,q}}<\infty\},$$ \end{defin}
\noindent where
$$\|f\|_{{B}^s_{p,q}}=\left\{\begin{array}{l}
\displaystyle\bigg(\sum_{j\ge 0}2^{jsq}\|\Delta_j f\|_p^q\bigg)^{\frac 1
q}+\|S_0(f)\|_p,\quad \hbox{for}\quad q<\infty,\\
\displaystyle\sup_{j\ge 0}2^{js}\|\Delta_jf\|_p+\|S_0(f)\|_p, \quad
\hbox{ for} \quad q=\infty.
\end{array}\right.
$$
If $s>0$, then ${B}^s_{p,q}=L^p\cap\dot{B}^s_{p,q}$ and $\|f\|_{B^s_{p,q}}\approx\|f\|_{p}+\|f\|_{\dot{B}^s_{p,q}}.$
We refer to \cite{B-L,Tr} for  details.

The following Definition \ref{def2.3} gives the mixed time-space
Besov space dependent on Littlewood-Paley decomposition (cf.
\cite{C-L}).

\begin{defin}\label{def2.3}
Let $u(t,x)\in \mc{S}'(\mr^{n+1})$, $s\in\mr$,$1\le p,\ q,\ \rho
\le\infty$.  We say that $u(t,x)\in
\widetilde{L}^\rho\Big(I;\dot{B}^s_{p,q}(\mr^n)\Big)$ if and only
if
$$
2^{js}\|\tr_ju\|_{L^\rho(I;L^p)}\in l^q,
$$
and we define
 \be\label{2.16}
 \|u\|_{\widetilde{L}^\rho(I;\dot{B}^s_{p,q})}\triangleq\bigg(\sum_{j\in\mb{Z}}2^{jsq}
 \|\tr_ju\|^q_{L^\rho(I;L^p)}\bigg)^{1/q}.
 \ee
\end{defin}

For the convenience we also recall the definition of Bony's
para-product formula which gives the decomposition of the product
 of two functions $f(x)$ and $g(x)$ (cf.
\cite{Bony,Cannone}).


\begin{defin}\label{def2.4}
The para-product of two functions $f$ and $g$ is defined by
\begin{eqnarray}
T_gf=\sum_{i\le j-2}\triangle_ig\triangle_jf=\sum_{j\in
\mathbb{Z}}S_{j-1}g\triangle_jf.
\end{eqnarray}
The remainder of the para-product is defined by
\begin{eqnarray}
R(f,g)=\sum_{|i-j|\le1}\triangle_ig\triangle_jf.
\end{eqnarray}
Then Bony's para-product formula reads
\begin{eqnarray}\label{2.18}
f\cdot g=T_gf+T_fg+R(f,g).
\end{eqnarray}
\end{defin}

Using  Bony's para-product formula and the definition of
homogeneous Besov space, one can prove the following trilinear
estimates, for details, see \cite{G-P}.

\begin{prop}\label{prop2.1}
Let $n\ge 2$ be the spatial dimension and let $r$ and $\sigma$ be
two real numbers such that $2\le r<\infty$, $2<\sigma<\infty$ and
$\frac nr+\frac2\sigma>1$. Define the trilinear form as
 \be
T(a,b,c)=\int^t_0\int_{\mr^n}(a(s, x)\cdot\nabla  b(s, x))\cdot c(s, x)\md
x\md s,
 \ee
 for $a,b\in L^\infty([0,\infty);L^2(\mr^n))\cap
 L^2([0,\infty);\dot{H}^1(\mr^n))$ and
 $c\in L^\sigma([0,T];\dot{B}^{n/r+2/\sigma-1}_{r,\sigma}(\mr^n))$, $0<t\le
 T$. Then $T(a,b,c)$ is continuous and satisfies estimates as
 follows:
  \be\nonumber
|T(a,b,c)|&\lesssim&
\|a\|^{1/\sigma}_{L^\infty(\mr^+;L^2)}\|\nabla
a\|^{1-1/\sigma}_{L^2(\mr^+;L^2)}\|b\|^{1/\sigma}_{L^\infty(\mr^+;L^2)}\|\nabla
b\|^{1-1/\sigma}_{L^2(\mr^+;L^2)}\|c\|_{L^\sigma([0,T];\dot{B}^{\frac
nr+\frac2\sigma-1}_{r,\sigma})}\\
\nonumber &\ & +\|\nabla
a\|_{L^2(\mr^+;L^2)}\|b\|^{2/\sigma}_{L^\infty(\mr^+;L^2)}\|\nabla
b\|^{1-2/\sigma}_{L^2(\mr^+;L^2)}\|c\|_{L^\sigma([0,T];\dot{B}^{\frac
nr+\frac2\sigma-1}_{r,\sigma})}\\ \label{2.21} &\ &+
\|a\|^{2/\sigma}_{L^\infty(\mr^+;L^2)}\|\nabla
a\|^{1-2/\sigma}_{L^2(\mr^+;L^2)}\|\nabla
b\|_{L^2(\mr^+;L^2)}\|c\|_{L^\sigma([0,T];\dot{B}^{\frac
nr+\frac2\sigma-1}_{r,\sigma})},
  \ee
and
 \be\nonumber
|T(a,b,c)|&\le& C(\ep)(\|\nabla a\|_{L^2(\mr^+;L^2)}+\|\nabla
b\|_{L^2(\mr^+;L^2)})\\ \label{2.22} &\ &
+C(\ep^{-1})\int^t_0(\|a(s)\|^2_{L^2}+\|b(s)\|^2_{L^2})
\|c(s)\|^\sigma_{\dot{B}^{\frac nr+\frac2\sigma-1}_{r,\sigma}}\md
s.
 \ee
 In particular,
 \be\label{2.23}
|T(a,a,c)|\le C(\ep)\|\nabla
a\|_{L^2(\mr^+;L^2)}+C(\ep^{-1})\int^t_0\|a(s)\|^2_{L^2}\|c(s)\|^\sigma_{\dot{B}^{\frac
nr+\frac2\sigma-1}_{r,\sigma}}\md s.
 \ee
Here $C(\ep)$ and $C(\ep^{-1})$ are constants that can be arranged
by $\ep$ and $\frac1\ep$, respectively, for $\ep>0$.
\end{prop}

\begin{rem}
In reference \cite{G-P} authors only proved the estimates
(\ref{2.21}) and (\ref{2.23}). Actually the proof also implies the
estimate (\ref{2.22}).
\end{rem}

Next we give  the time-space estimate of the  heat semigroup
$u(t,x)=S(t)u_0\triangleq \mx{e}^{-t\tr}u_0(x)$, which has been
proved  in  \cite{G-P}.  But the proof has a misprint that is the
inequality (4.3) in \cite{G-P} should be
 \be
\|u\|_{L_t^{p^\mp}(I; \dot{B}^{s\pm}_{2,2})}\lesssim
\|u\|^{1-\frac2{p\mp}}_{\widetilde{L}_t^\infty(I; \dot{B}^{\pm\ep}_{2,2})}
\||\nabla|^{1\pm\ep}u\|^{\frac2{p^\mp}}_{L^2_{t,x}(I\times\mathbb{R}^n)},
 \ee
where $I\subset [0, \infty)$ or $I=[0,\infty)$.

\begin{prop} \label{prop2.2}
Let $2<p<\infty$, $u_0(x)\in L^2(\mr^n)$. Denote
$u(t,x)=S(t)u_0(x)$, then we have
 \be
\|u\|_{L^{p,2}_t(I; L^q_x)}\le C\|u_0\|_{L^2},
 \ee
for $\frac2p+\frac nq=\frac n2$, $L^{p,2}_t(I)$ denotes Lorentz space with respect to $t\in I$.
\end{prop}

The following propositions describe the H\"older's and Young's
inequalities in Lorentz spaces, which will be used in this paper,
for their proofs we  refer to \cite{O'neil}.

\begin{prop}\label{mtholder}
(\textbf{Generalized H\"{o}lder's inequality}) Let $1< p_1,\ p_2,\
r<\infty$, such that
\begin{eqnarray*}
\frac1r=\frac 1{p_1}+\frac 1{p_2}<1,
\end{eqnarray*}
 and $1\le q_1,\ q_2,\ s\le\infty$ with
\begin{eqnarray*}
\frac1{q_1}+\frac1{q_2}\ge\frac1s.
\end{eqnarray*}
If  $f\in L^{p_1,\:q_1},\ g\in L^{p_2,\:q_2}$, then $h=fg\in
L^{r,\:s}$  such that
\begin{eqnarray}\label{holder}
\|h\|_{(r,\:s)}\le r'\|f\|_{(p_1,\:q_1)}\|g\|_{(p_2,\:q_2)},
\end{eqnarray}\\
where $r'$ stands for the dual to $r$, i.e.
$\frac1r+\frac1{r'}=1$.
\end{prop}

\begin{prop}\label{mtyoung}
(\textbf{Generalized Young's inequality})

Let $1< p_1,\ p_2,\  r<\infty$ such that
\begin{eqnarray*}
\frac1{p_1}+\frac1{p_2}>1,\ \frac1r=\frac 1{p_1}+\frac 1{p_2}-1,
\end{eqnarray*}
and $1\le q_1,\ q_2,\ s\le\infty$ with
\begin{eqnarray*}
\frac1{q_1}+\frac1{q_2}\geq\frac1s.
\end{eqnarray*}
If  $f\in
L^{p_1,\:q_1},\ g\in L^{p_2,\:q_2}$, then $h=f*g\in L^{r,\:s}$
with
\begin{eqnarray}\label{young}
\|h\|_{(r,\:s)}\le3r\|f\|_{(p_1,\:q_1)}\|g\|_{(p_2,\:q_2)}.
\end{eqnarray}
In particular, we have the weak Young's inequality
\begin{eqnarray}\label{young1}
\|h\|_{(r,\infty)}\le C(p,q)\|f\|_{(p,\infty)}\|g\|_{(q,\infty)},
\end{eqnarray}
where $1<p,\ q,\ r<\infty$ and $\frac1r=\frac1p+\frac1q-1$.
\end{prop}

\begin{prop}\label{Conv}
Let  $1\le q_1\le\infty$ and
$1\le q_2\le\infty$ satisfy $\frac1{q_1}+\frac1{q_2}\ge 1$, $p$
and $p'$ be conjugate indices, i.e. $\frac1p+\frac1{p'}=1$. If  $f(x)\in
L^{p,\; q_1}$ and $g(x)\in L^{p',\; q_2}$, then
$h(x)=f*g\in L^\infty$ such that
 \be\label{ConvIneq}
\|h\|_\infty\le
\|f\|_{(p,\:q_1)}\|g\|_{(p',\:q_2)}.
 \ee
\end{prop}


\section{Well-posedness in Besov spaces: Case $n\ge 2$}
\setcounter{equation}{0}

 This section is devoted to the proof of Theorem \ref{thm1.1}.
One easy sees that (\ref{1.1})-(\ref{1.5}) can be rewritten as
 \be\label{3.1}
u_t-\tr u+\mb{P}\nabla\cdot(u\otimes u)-\mb{P}\nabla\cdot(b\otimes
b)=0,\\ \label{3.2}
b_t-\tr b+\mb{P}\nabla\cdot(u\otimes
b)-\mb{P}\nabla\cdot(b\otimes
b)=0,\\
\mx{\rm div} u=\mx{\rm div} b=0,\\
\label{3.4} u(0,x)=u_0(x),\ b(0,x)=b_0(x).
 \ee
or their integral form

 \be\label{3.1add}
u= e^{t\tr}u_0-\int_0^t e^{(t-s)\tr}\Big[\mb{P}\nabla\cdot(u\otimes u)-\mb{P}\nabla\cdot(b\otimes
b)\Big]ds,\\
\label{3.2add}
b= e^{t\tr}b_0-\int_0^te^{(t-s)\tr} \Big[\mb{P}\nabla\cdot(u\otimes
b)-\mb{P}\nabla\cdot(b\otimes
u)\Big]ds.
 \ee
Here $\mb{P}$ stands for the Leray projector onto divergence free
vector field.

\subsection{Linear and nonlinear estimates}

To prove the results of global or local well-posedness of the
Cauchy problem (\ref{3.1})-(\ref{3.4}) or (\ref{3.1add})-(\ref{3.2add}) in Besov space
$\dot{B}^{n/p-1}_{p,r}(\mr^n)$,  we need to establish linear and nonlinear
estimates in framework of mixed space-time space by Fourier localization.
First we consider the solution to linear parabolic
equation

 \be\label{3.7}
 \begin{cases}
u_t-\tr u=f(t,x),\\
u(0,x)=u_0.
 \end{cases}
 \ee
Applying frequency projection operator $\tr_j$ to  both sides of
(\ref{3.7}), one arrives at
 \be\label{3.8}
\frac{\pt}{\pt t}(\tr_ju)+\tr(\tr_ju)=\tr_j f.
 \ee
Multiplying $|\tr_ju|^{p-2}\tr_ju$ on both sides of
(\ref{3.8}),  we obtain
  \be\label{3.9}
  \frac {\pt}{\pt t}\tr_j u|\tr_ju|^{p-2}\tr_ju-\tr\tr_j u|\tr
  u|^{p-2}\tr_ju=\tr_jf|\tr_ju|^{p-2}\tr_ju.
 \ee
We integrate both sides of (\ref{3.9}) and apply the divergence
theorem to obtain
 \be
 \frac1p\frac{\md}{\md t}\|\tr_j
 u\|^p_p+\int_{\mr^n}\nabla\tr_ju\cdot\nabla(|\tr_ju|^{p-2}\tr_ju)\md
 x\le \|\tr_jf\|_p\|\tr_ju\|^{p-1}_p.
\ee
 Since
 \be\nonumber
 \int_{\mr^n}\nabla\tr_ju\cdot\nabla(|\tr_ju|^{p-2}\tr_ju)\md x &=& (p-1)\int_{\mr^n}
 |\tr_ju|^{p-2}|\nabla\tr_ju|^2\md x
\\ \nonumber &=&\frac {4(p-1)}{p^2}\int_{\mr^n}|\nabla(|\tr_ju|^{\frac
p2})|^2\md x=\|\nabla(|\tr_ju|^{\frac p2})\|^2_2\\
 &\ge&c_p2^{2j}\|\tr_ju\|^p_p.
\ee
 We have
\be \label{3.10}
 \frac{\md}{\md t}\|\tr_ju\|_p+2^{2j}c_p\|\tr_ju\|_p\le
 \|\tr_jf\|_p.
\ee
 Integrating both sides of (\ref{3.10}) with respect to $t$ we
 arrive at
\be\label{3.11}
 \|\tr_ju\|_p\le
 \mx{e}^{-c_p2^{2j}t}\|\tr_ju_0(x)\|_p+\mx{e}^{-2^{2j}c_pt}*(\|\tr_jf\|_p\chi(\tau)),
\ee
 where $\chi(\tau)$ is a character function
\be\label{3.12}
 \chi(\tau)=
\begin{cases}
1,\ &\mx{if }0\le \tau\le t,\\
0,\ &\mx{if }\mx{others}.
\end{cases}
\ee
 Taking $L^q-$norm with respect to $t$ in interval $I$ in both
 sides of (\ref{3.12}), by Young inequality one has
\be\label{3.14}
 \|\tr_ju\|_{L^q(I;L^p)}\le
 c_p^{-\frac1q}2^{-\frac{2j}q}\|\tr_ju_0(x)\|_p+
 C(p,q)2^{-\frac{2j}{q'}}\|\tr_jf\|_{L^{q/2}(I;L^p)}.
\ee
 Here $\frac1q+\frac1{q'}=1$.
 Multiplying $2^{js+\frac{2j}q}$ on both sides of (\ref{3.14}) and
 taking $l^r-$norm with respect to $j$ yields
\be
 \|u\|_{\widetilde{L}^q(I;\dot{B}^{s+2/q}_{p,r})}\le
 C(p,q)\bigg(\|u_0\|_{\dot{B}^s_{p,r}}+
 \|f\|_{\widetilde{L}^{q/2}(I;\dot{B}^{s+4/q-2}_{p,r})}\bigg).
\ee
Thus we arrive at

\begin{lem}\label{lem3.1}
Let $1\le p<\infty$, $2\le q\le \infty$, $1\le r\le \infty$ and
$s\in \mr$. Assume $u(t,x)$ is a solution to the Cauchy problem
(\ref{3.7}). Then there exists a constant $C$ depending on $p,\
q,\ n$ so that
 \be
 \|u\|_{\widetilde{L}^q(I;\dot{B}^{s+2/q}_{p,r})}\le
 C\bigg(\|u_0\|_{\dot{B}^s_{p,r}}+
 \|f\|_{\widetilde{L}^{q/2}(I;\dot{B}^{s+4/q-2}_{p,r})}\bigg).
\ee
 In particular, if $\frac q2\le r\le q$  we have
\be
 \|u\|_{L^q(I;\dot{B}^{s+2/q}_{p,r})}\le
 C\bigg(\|u_0\|_{\dot{B}^s_{p,r}}+
 \|f\|_{L^{q/2}(I;\dot{B}^{s+4/q-2}_{p,r})}\bigg),
\ee
 by Minkowski inequality.
 \end{lem}

\begin{rem} For $s\in \mathbb{R}$ and
$1\le p,r\le\infty$, $ ( {\dot B}^s_{p,r}, \|\cdot\|_{{\dot B}^s_{p,r}})$ is a normed space. It is easy to check that
$( {\dot B}^s_{p,r}, \|\cdot\|_{{\dot B}^s_{p,r}})$   is a  Banach space if and only if
$s<\frac{n}p$ or $s=\frac{n}p$, $r=1$.
\end{rem}

Using Bony's para-product decomposition we study the bilinear estimates.
 Consider two tempered distributions $u(t,x)$ and
$v(t,x)$, then
 \be
uv=T_uv+T_vu+R(u,v).
 \ee
First we deal with the para-product term $T_uv$ or $T_vu$ as
following lemma.

\begin{lem}\label{lem3.2}
(1) Let  $\dot{B}^s_{p,r}(\mr^n)$ be a
Banach space, then
 \be
  \|T_uv\|_{\widetilde{L}^{q/2}(I;\dot{B}^s_{p,r})}\le
  \|u\|_{L^q(I;L^\infty)}\|v\|_{\widetilde{L}^q(I;\dot{B}^s_{p,r})}.
  \ee

(2) Let $s_1<0$ and $\frac1r=\frac1{r_1}+\frac1{r_2}$, and
 $\dot{B}^{s_2}_{p,r_2}(\mr^n)$
be a Banach space. Then
 \be
\|T_uv\|_{\widetilde{L}^{q/2}(I;\dot{B}^{s_1+s_2}_{p,r})}\le
 C \|u\|_{\widetilde{L}^q(I;\dot{B}^{s_1}_{\infty,r_1})}
  \|v\|_{\widetilde{L}^q(I;\dot{B}^{s_2}_{p,r_2})}.
 \ee
\end{lem}

\begin{proof}
 By the definition of $\widetilde{L}^q(I;\dot{B}^{s}_{p,r})$
and H\"older inequality, direct computation yields
 \be\nonumber
\|T_uv\|_{\widetilde{L}^{q/2}(I;\dot{B}^s_{p,r})}&=&
\bigg(\sum_{j\in\mb{Z}}2^{jsr}\|S_{j-1}u\tr_jv\|^r_{L^{q/2}(I;L^p)}\bigg)^{1/r}\\\nonumber
&\le&
\|u\|_{L^q(I;L^\infty)}\bigg(\sum_{j\in\mb{Z}}2^{jsr}\|\tr_jv\|^r_{L^q(I;L^p)}\bigg)^{1/r}
\\ &\le& \|u\|_{L^q(I;L^\infty)}\|v\|_{\widetilde{L}^q(I;\dot{B}^s_{p,r})}.
 \ee

\noindent  Noting that the equivalent definition of negative index Besov space
 \be
\bigg(\sum_{j\in\mb{Z}}2^{jsr}\|S_ju\|^r_{L^q(I;L^p)}\bigg)^{1/r}\simeq
\bigg(\sum_{j\in\mb{Z}}2^{jsr}\|\tr_ju\|^r_{L^q(I;L^p)}\bigg)^{1/r}
 \ee
for $s<0$ (cf. \cite{Cannone}),  we can
derive similarly by  H\"older inequality
 \be\nonumber
\|T_uv\|_{\widetilde{L}^{q/2}(I;\dot{B}^{s_1+s_2}_{p,r})}&\le&
\bigg(\sum_{j\in\mb{Z}}2^{j(s_1+s_2)r}\|S_{j-1}u\|^r_{L^q(I;L^\infty)}
\|\tr_jv\|^r_{L^q(I;L^p)}\bigg)^{1/r}\\ &\le& C
\|u\|_{\widetilde{L}^q(I;\dot{B}^{s_1}_{\infty,r_1})}
  \|v\|_{\widetilde{L}^q(I;\dot{B}^{s_2}_{p,r_2})}.
 \ee
\end{proof}

Next we estimate the remainder of para-product decomposition.

\begin{lem}\label{lem3.3}
Let $s_1,\ s_2\in \mr$,  $1\le p_1,\ p_2,\ p,\ r_1,\ r_2,\
r\le\infty$ and $2\le q\le\infty$ such that
 \be
\frac1p=\frac1{p_1}+\frac1{p_2},\quad
\frac1r=\frac1{r_1}+\frac1{r_2}
 \ee
and $\widetilde{L}^q(I;\dot{B}^{s_1}_{p_1,r_1})$,
$\widetilde{L}^q(I;\dot{B}^{s_2}_{p_2,r_2})$ and
$\widetilde{L}^{q/2}(I;\dot{B}^{s_1+s_2}_{p,r})$ are Banach
spaces. Assume $0<s_1+s_2<\frac np$, then
 \be
\|R(u,v)\|_{\widetilde{L}^{q/2}(I;\dot{B}^{s_1+s_2}_{p,r})}\le C
\|u\|_{\widetilde{L}^q(I;\dot{B}^{s_1}_{p_1,r_1})}
\|v\|_{\widetilde{L}^q(I;\dot{B}^{s_2}_{p_2,r_2})}.
 \ee
Moreover, if $s_1+s_2=0$ and $\frac1{r_1}+\frac1{r_2}=1$, then one
has
 \be
\|R(u,v)\|_{\widetilde{L}^{q/2}(I;\dot{B}^0_{p,\infty})}\le C
\|u\|_{\widetilde{L}^q(I;\dot{B}^{s_1}_{p_1,r_1})}
\|v\|_{\widetilde{L}^q(I;\dot{B}^{s_2}_{p_2,r_2})}.
 \ee
 If $s_1+s_2=\frac np$ and $r=1$, then
 \be
\|R(u,v)\|_{\widetilde{L}^{q/2}(I;\dot{B}^{n/p}_{p,1})}\le C
\|u\|_{\widetilde{L}^q(I;\dot{B}^{s_1}_{p_1,r_1})}
\|v\|_{\widetilde{L}^q(I;\dot{B}^{s_2}_{p_2,r_2})}.
 \ee
\end{lem}

\begin{proof}
Write
 \be
R(u,v)=\sum_{j'\in\mb{Z}}\sum_{k=-1}^1\tr_{j'}u\tr_{k+j'}v\triangleq
\sum_{j'\in\mb{Z}}R_{j'}.
 \ee
Since supp$[\mc{F}(\tr_{j'}u\tr_{k+j'}v)]\subseteqq
\{|\xi|\le\frac83 2^{j'}(1+2^k)\}$ and
supp$[\mc{F}(\tr_jf)]\subseteqq\{\frac34 2^j\le |\xi|\le \frac83
2^j\}$, it follows that
 \be
\tr_jR(u,v)=\sum_{j'\ge j-4}\tr_jR_{j'}.
 \ee
A straightforward calculation shows that
 \be\nonumber
&\ & 2^{j(s_1+s_2)}\|\tr_jR_{j'}\|_{L^{q/2}(I;L^p)}\le
2^{j(s_1+s_2)}\sum_{k=-1}^1\|\tr_{j'+k}u\|_{L^q(I;L^{p_1})}\|\tr_{j'}v\|_{L^q(I;L^{p_2})}\\
\label{3.31} &\le& \sum_{k=-1}^1
2^{-(j'-j)(s_1+s_2)}2^{j's_1}\|\tr_{j'+k}u\|_{L^q(I;L^{p_1})}2^{j's_2}
\|\tr_{j'}v\|_{L^q(I;L^{p_2})}.
 \ee
Using the estimate (\ref{3.31}) and the definition (\ref{2.16}) of
mixed time-space Besov space one has
 \be\label{3.32}
&\ &\|R(u,v)\|_{\widetilde{L}^{q/2}(I;\dot{B}^{s_1+s_2}_{p,r})}\\
\nonumber &\le& \bigg(\sum_{j\in\mb{Z}}\Big(\sum_{l\le
4}\sum^1_{k=-1}2^{l(s_1+s_2)-s_1k}2^{s_1(j+k-l)}\|\tr_{j+k-l}u\|_{L^qL^{p_1}}
2^{s_2(j-l)}\|\tr_{j-l}v\|_{L^qL^{p_2}}\Big)^r\bigg)^{1/r}
 \ee
In view of Minkowski and H\"older inequalities,  we obtain
 \be\nonumber
\|R(u,v)\|_{\widetilde{L}^{q/2}(I;\dot{B}^{s_1+s_2}_{p,r})}&\le&
\sum_{k=-1}^1 2^{-s_1k}\sum_{l\le
4}2^{l(s_1+s_2)}\|u\|_{\widetilde{L}^q(I;\dot{B}^{s_1}_{p_1,r_1})}
\|v\|_{\widetilde{L}^q(I;\dot{B}^{s_2}_{p_2,r_2})}\\ &\le&
C\|u\|_{\widetilde{L}^q(I;\dot{B}^{s_1}_{p_1,r_1})}
\|v\|_{\widetilde{L}^q(I;\dot{B}^{s_2}_{p_2,r_2})}.
 \ee
In particular, if $s_1+s_2=0$, we first apply Minkowski
inequality then  H\"older inequality  to the right of (\ref{3.31}),
it follows that
 \be\nonumber
\|R(u,v)\|_{\widetilde{L}^{q/2}(I;\dot{B}^0_{p,\infty})}&\le&
\sum_{k=-1}^1 \|u\|_{\widetilde{L}^q(I;\dot{B}^{s_1}_{p_1,r_1})}
\|v\|_{\widetilde{L}^q(I;\dot{B}^{s_2}_{p_2,r_2})}\\ &\le&
C\|u\|_{\widetilde{L}^q(I;\dot{B}^{s_1}_{p_1,r_1})}
\|v\|_{\widetilde{L}^q(I;\dot{B}^{s_2}_{p_2,r_2})}.
 \ee
  The proof of Lemma \ref{lem3.3} is thus complete.
\end{proof}

By means of the fact
$L^\infty(\mr^n)\hookrightarrow
\dot{B}^0_{\infty,\infty}(\mr^n)$ and
$$\widetilde{L}^q(I;\dot{B}^{n/p}_{p,1}(\mr^n))\hookrightarrow
L^q(I;\dot{B}^{n/p}_{p,1})\hookrightarrow L^q(I;L^\infty(\mr^n)), \quad
q\ge 1,$$
we apply  Lemma \ref{lem3.3} with $p_1=r_1=\infty$, $s_1=0$ and Lemma \ref{lem3.2} to get
the following  results.

\begin{cor}\label{cor3.1}
Let $s$ be a real number such that $s<\frac np$, $q\ge 2$ and
$1\le p,\ r\le \infty$, one has
 \be
\|uv\|_{\widetilde{L}^{q/2}(I;\dot{B}^s_{p,r})}\le
C\|u\|_{L^q(I;L^\infty)}\|v\|_{\widetilde{L}^q(I;\dot{B}^s_{p,r})}
+\|u\|_{\widetilde{L}^q(I;\dot{B}^s_{p,r})}\|v\|_{L^q(I;L^\infty)}
 \ee
 and
 \be
\|uv\|_{\widetilde{L}^{q/2}(I;\dot{B}^{n/p}_{p,1})}\le
C\|u\|_{\widetilde{L}^q(I;\dot{B}^{n/p}_{p,1})}
\|v\|_{\widetilde{L}^q(I;\dot{B}^{n/p}_{p,1})}.
 \ee
\end{cor}

\begin{cor}\label{cor3.2}
Let $s_1,\ s_2\in \mr$, $1\le p_k,\ r_k\le\infty$ and $1\le p,\
r\le\infty$ such that
 \be
s_k<\frac n{p_k},\ \frac1{r_1}+\frac1{r_2}=\frac1r,\
p\ge\max(p_1,\ p_2)
 \ee
for $k=1,\ 2$. If
$s_1+s_2>n\Big(\frac1{p_1}+\frac1p_2-\frac1p\Big)$, then
 \be\label{3.38}
\|uv\|_{\widetilde{L}^{q/2}(I;\dot{B}^{s_1+s_2-n(\frac1{p_1}+\frac1p_2-\frac1p)}_{p,r})}
\le C\|u\|_{\widetilde{L}^q(I;\dot{B}^{s_1}_{p_1,r_1})}
\|v\|_{\widetilde{L}^q(I;\dot{B}^{s_2}_{p_2,r_2})}.
 \ee
\end{cor}

\begin{proof}
Let $\widetilde{s_1}=s_1-\frac n{p_1}$ and
$\widetilde{s_2}=s_2-\frac n{p_2}+\frac np$, then
$\widetilde{s_1}<0$. Applying Lemma \ref{lem3.2} to para-product
$T_uv$ and $T_vu$, Lemma \ref{lem3.3} to the remainder $R(u,v)$ we
obtain that
 \be
\|uv\|_{\widetilde{L}^{q/2}(I;\dot{B}^{\widetilde{s_1}+\widetilde{s_2}}_{p,r})}\le
C\|u\|_{\widetilde{L}^q(I;\dot{B}^{\widetilde{s_1}}_{\infty,r_1})}
\|v\|_{\widetilde{L}^q(I;\dot{B}^{\widetilde{s_2}}_{p,r_2})}.
 \ee
Noting the embedding relations
 \be
\widetilde{L}^q(I;\dot{B}^{s_1}_{p_1,r_1}(\mr^n))
\hookrightarrow\widetilde{L}^q(I;\dot{B}^{\widetilde{s_1}}_{\infty,r_1}(\mr^n)),\
\widetilde{L}^q(I;\dot{B}^{s_2}_{p_2,r_2}(\mr^n))\hookrightarrow
\widetilde{L}^q(I;\dot{B}^{\widetilde{s_2}}_{p,r_2}(\mr^n)),
 \ee
we complete the proof of the estimate (\ref{3.38}).
\end{proof}

\subsection{Well-posedness in Besov spaces and uniqueness of weak and strong solutions}

In this subsection we first prove Theorem \ref{thm1.1}, i.e.  the small global well-posedness of
the Cauchy problem (\ref{3.1})-(\ref{3.4}) or (\ref{3.1add})-(\ref{3.2add}) and local
well-posedness in Besov space $\dot{B}^{n/p-1}_{p,r}(\mr^n)$ for $n\ge 2$.
Then we establish the stability
result of the Leray weak
solution and strong solution which implies  Theorem\ref{1.2}.

{\bf The proof of Theorem\ref{thm1.1}. }
Without loss of generality we can prove the case $r\not=\infty$.
By Lemma \ref{lem3.1} we first prove the following bilinear
estimate.
 \be
\|\mb{P}\nabla\cdot(u\otimes
b)\|_{\widetilde{L}^{q/2}(I;\dot{B}^{s+4/q-2}_{p,r})}\le
C\|u\|_{\widetilde{L}^q(I;\dot{B}^{s+2/q}_{p,r})}
\|b\|_{\widetilde{L}^q(I;\dot{B}^{s+2/q}_{p,r})}.
 \ee
Indeed, by the boundedness of Calder\'{o}n-Zygmund singular
integral operator on the space
$\widetilde{L}^{q/2}(I;\dot{B}^{s+4/q-2}_{p,r})$ and the Sobolev
embedding theorem one has
 \be
\|\mb{P}\nabla\cdot(u\otimes
b)\|_{\widetilde{L}^{q/2}(I;\dot{B}^{s+4/q-2}_{p,r})}\le
C\|u\otimes
b\|_{\widetilde{L}^{q/2}(I;\dot{B}^{s+4/q-1}_{p,r/2})}.
 \ee
Taking $s_1=s_2=s_p+\frac2q$ and
$r_1=r_2=\frac r2$ in Corollary \ref{cor3.2} we have
 \be\label{3.43}
\|u\otimes
b\|_{\widetilde{L}^{q/2}(I;\dot{B}^{s+4/q-1}_{p,r/2})}\le C
\|u\|_{\widetilde{L}^q(I;\dot{B}^{s+2/q}_{p,r})}
\|b\|_{\widetilde{L}^q(I;\dot{B}^{s+2/q}_{p,r})}.
 \ee
Thus we get estimates of solution to equations
(\ref{3.1add})-(\ref{3.2add}) as follows
 \be\label{3.44}
 \begin{cases}
\|u\|_{\widetilde{L}^q(I;\dot{B}^{s+2/q}_{p,r})}\le
C\|u_0\|_{\dot{B}^{s+2/q}_{p,r})}+C\|b\|^2_{\widetilde{L}^q(I;\dot{B}^{s+2/q}_{p,r})}
+C\|b\|^2_{\widetilde{L}^q(I;\dot{B}^{s+2/q}_{p,r})},\\
\|b\|_{\widetilde{L}^q(I;\dot{B}^{s+2/q}_{p,r})}\le
C\|b_0\|_{\dot{B}^{s+2/q}_{p,r})}+2C\|u\|_{\widetilde{L}^q(I;\dot{B}^{s+2/q}_{p,r})}
\|b\|_{\widetilde{L}^q(I;\dot{B}^{s+2/q}_{p,r})}.
\end{cases}
 \ee
For convenience we write
 \be
\|(u,b)\|_{\widetilde{L}^q(I;\dot{B}^{s+2/q}_{p,r})}=
\|u\|_{\widetilde{L}^q(I;\dot{B}^{s+2/q}_{p,r})}
+\|b\|_{\widetilde{L}^q(I;\dot{B}^{s+2/q}_{p,r})},
 \ee
thus estimate (\ref{3.44}) can be written consequently  as
 \be\label{3.46}
\|(u,b)\|_{\widetilde{L}^q(I;\dot{B}^{s+2/q}_{p,r})}\le
C\|(u_0,b_0)\|_{\dot{B}^{s+2/q}_{p,r})}+C\|(u,b)\|^2
_{\widetilde{L}^q(I;\dot{B}^{s+2/q}_{p,r})}.
 \ee
Let $(u_1,b_1)$ and $(u_2,b_2)$ be two solutions to the Cauchy
problem (\ref{3.1add})-(\ref{3.2add}) with the same data $(u_0,b_0)$.
Arguing  similarly as in deriving (\ref{3.46}), one has
 \be
\|(u_1-u_2,b_1-b_2)\|_{\widetilde{L}^q(I;\dot{B}^{s+2/q}_{p,r})}\le&
C (\|(u_1,b_1)\|_{\widetilde{L}^q(I;\dot{B}^{s+2/q}_{p,r})}+
\|(u_2,b_2)\|_{\widetilde{L}^q(I;\dot{B}^{s+2/q}_{p,r})})\nonumber\\
 &\times
\|(u_1-u_2,b_1-b_2)\|_{\widetilde{L}^q(I;\dot{B}^{s+2/q}_{p,r})}.\label{3.46add}
 \ee
We finally, apply the Banach contraction mapping principle to
nonlinear operator ${\cal T}$ defined by the right sides of (\ref{3.1add})-(\ref{3.2add}) in
 a closed set $E$
 \be
 E=\Big\{(f,g):\;  \|(f, g)\|_{\widetilde{L}^q(I;\dot{B}^{s+2/q}_{p,r})}\le 2K_0\Big\},
 \ee
where $K_0$ is the constant dependent on the local existence time $T$
in local existence case or on the initial datum norm
$\|(u_0,b_0)\|_{\dot{B}^{n/p-1}_{p,r}}\ll1$ in global case. In face, we can
choose the existence time $T$ small enough (cf. estimate
(\ref{3.11})) or the norm $\|(u_0,b_0)\|_{\dot{B}^{n/p-1}_{p,r}}$
small enough so that ${\cal T}$ is a contraction mapping on $E$ by (\ref{3.46}) and (\ref{3.46add}).
  Thus  a standard argument together with Remark \ref{rem3.2}  shows Theorem \ref{1.1} and some further regularity
  of solution $(u,b)$.

\begin{rem}\label{rem3.2}
(1) By Young inequality we have
 \be
\|S(t)(u_0,b_0)\|_{\widetilde{L}^\infty(I;\dot{B}^{n/p-1}_{p,r})}\le
\|(u_0,b_0)\|_{\dot{B}^{n/p-1}_{p,r}}.
 \ee
Consequently, combining the bilinear estimate (\ref{3.43}) it
follows that
 \be
(u,b)\in
\widetilde{L}^\infty(I;\dot{B}^{n/p-1}_{p,r}(\mr^n))\hookrightarrow
L^\infty(I;\dot{B}^{n/p-1}_{p,r}(\mr^n)).
 \ee

 (2) If $p>n$ and $1\le r\le \infty$, using the equivalent
characterization of negative index homogeneous Besov space
\cite{M-Y-Z}, we have
 \be\nonumber
\sup_{t>0}t^{\frac12-\frac n{2p}+\frac \alpha2}\|\nabla^\alpha
S(t)(u_0,b_0)\|_{L^p}&\le&C\|\nabla^\alpha(u_0,b_0)\|_{\dot{B}^{n/p-1-\alpha}_{p,\infty}}\\
&\le&C\|(u_0,b_0)\|_{\dot{B}^{n/p-1}_{p,r}},
 \ee
for $\alpha=0,\ 1$. Denote that
 \be
B(u,v)=\int^t_0S(t-s)\mb{P}\nabla\cdot(u\otimes v)\md s
 \ee
 A straightforward calculation shows that
 \be\nonumber
\sup_{t>0}t^{\frac12-\frac n{2p}}\|B(u,v)\|_{L^p}&\le&
C\sup_{t>0}t^{\frac12-\frac n{2p}}\int^t_0(t-s)^{-\frac12-\frac
n{2p}}\|u\|_{L^p}\|v\|_{L^p}\md s\\
&\le& C\sup_{t>0}t^{\frac12-\frac
n{2p}}\|u\|_{L^p}\sup_{t>0}t^{\frac12-\frac n{2p}}\|v\|_{L^p},
 \ee
for $\alpha=0$. Similarly, if $\alpha=1$, one has
 \be\nonumber
\sup_{t>0}t^{1-\frac n{2p}}\|\nabla B(u,v)\|_{L^p}&\le&
C\sup_{t>0}t^{1-\frac n{2p}}\int^t_0(t-s)^{-1-\frac
n{2p}}\|u\|_{L^p}\|v\|_{L^p}\md s\\
&\le& C\sup_{t>0}t^{\frac12-\frac
n{2p}}\|u\|_{L^p}\sup_{t>0}t^{\frac12-\frac n{2p}}\|v\|_{L^p}.
 \ee
Arguing similarly as in deriving Theorem \ref{thm1.1} we can prove
that the global solution $(u,b)$ to the Cauchy problem
(\ref{3.1})-(\ref{3.4}) for small datum enjoys the following
estimate
 \be\label{globalE1}
\|\nabla^{\alpha}(u,b)\|_{{\cal C}_{\frac12-\frac{n}{2p}}(\mathbb{R}^+; L^p)}\triangleq
\sup_{t>0}t^{\frac12-\frac n{2p}+\frac
\alpha2}\|\nabla^\alpha(u,b)\|_{L^p}\le C
\|(u_0,b_0)\|_{\dot{B}^{n/p-1}_{p,r}}\le \ep_0,
 \ee
for $\alpha=0$, or $1$, where $\ep_0>0$ is a small constant.
\end{rem}

To prove Theorem \ref{thm1.2},
we establish Proposition \ref{pro3.2} which describes one stability
 result for the weak and strong solutions
under some suitable conditions. As a direct consequence
we get the proof of Theorem \ref{thm1.2}, i.e. weak-strong uniqueness.

\begin{prop}\label{pro3.2}
Let $(u_0,b_0)$ and $(w_0,h_0)$ be the divergence free vector field in
$L^2(\mr^n)$, and $(w_0,h_0)$ be also in
$\dot{B}^{n/p-1}_{p,r}(\mr^n)$. Here assume $1\le p< \infty$ and
$2<r<\infty$ such that $\frac n{2p}+\frac2r>1$. Let $(w,g)\in
L^\infty(\mr^+;L^2(\mr^n))\cap L^2(\mr^+;\dot{H}^1(\mr^n))$ be
Leray weak solution associated with $(w_0,h_0)$, let $(u,b)$ be
the unique solution associated with $(u_0,b_0)$ with
$$(u,b)\in
L^r([0,T];\dot{B}^{\frac np+\frac2r-1}_{p,r}(\mr^n))\cap
L^\infty([0,T];L^2(\mr^n))\cap L^2([0,T];\dot{H}^1(\mr^n)), \quad \mbox{for some}\; T>0. $$
  Denote $(v,g)=(u,b)-(w,h)$, then we have for $0<t<T$
 \be\nonumber
\|(v,g)\|^2_{L^2}+\int^t_0\|\nabla(v(s),g(s))\|^2_{L^2}\md s&\le&
\exp\bigg(C\int^t_0\|(u(s),b(s))\|^r_{\dot{B}^{\frac
np+\frac2r-1}_{p,r}}\md s\bigg)\\
&\ & \times \|(u_0,b_0)-(w_0,h_0)\|^2_{L^2}.\label{3.55}
 \ee
\end{prop}

{\bf Proof of Proposition \ref{pro3.2}}. \quad To simplify the notation,  we write
 \be
(v_0,g_0)&\triangleq &(u_0,b_0)-(w_0,h_0),\\
 \|(\nabla v(s),\nabla
g(s))\|^2_{L^2}&\triangleq &\|\nabla v(s)\|^2_{L^2}+\|\nabla g(s)\|^2_{L^2},\\
\|(u(s),b(s))\|^r_{\dot{B}^{\frac
np+\frac2r-1}_{p,r}}&\triangleq&\|u(s)\|^r_{\dot{B}^{\frac
np+\frac2r-1}_{p,r}}+\|b(s)\|^r_{\dot{B}^{\frac
np+\frac2r-1}_{p,r}}.
 \ee
Subtracting the two equations satisfied by $(u,b)$ and $(w,h)$,
respectively, one has
 \be\label{3.56}
\pt_tv-\tr v+(v\cdot
\nabla)u+(w\cdot\nabla)v-[(g\cdot\nabla)b+(h\cdot\nabla)g]-\nabla(p_1-p_2)=0,\\\label{3.57}
\pt_tg-\tr
g+(v\cdot\nabla)b+(w\cdot\nabla)g-[(g\cdot\nabla)u+(h\cdot\nabla)v]=0.
 \ee
Formally, we may multiply equation (\ref{3.56}) and (\ref{3.57})
by $v$ and $g$ and integrate with respect to $x$ on $\mr^n$,
respectively, it follows
 \be\label{3.58}
\frac12\pt_t(v,v)+(\nabla v,\nabla v)+(v\cdot\nabla
u,v)+(w\cdot\nabla v,v)-\big[(g\cdot\nabla b,v)]+(h\cdot\nabla
g,v)\big]=0,\\ \label{3.59} \frac12\pt_t(g,g)+(\nabla g,\nabla
g)+(v\cdot\nabla b,g)+(w\cdot\nabla g,g)-\big[(g\cdot\nabla
u,g)]+(h\cdot\nabla v,g)\big]=0.
 \ee
Here $(\cdot,\cdot)$ denote the $L^2$-inner product. Integrating
(\ref{3.58}) and (\ref{3.59}), respectively, then summing up them
we arrive at
 \be\nonumber
\|(v,g)\|^2_{L^2}&+&2\int^t_0\|\nabla(v,g)\|^2_{L^2}\md s\le
\|(v_0,g_0)\|^2_{L^2}+2\int^t_0\int_{\mr^n}((v\cdot\nabla) vu\md x\md s\\
 &\ & +2\int^t_0\int_{\mr^n}((v\cdot\nabla)
gb-(g\cdot\nabla) vb-(g\cdot\nabla gu))\md x\md s,\label{3.64}
 \ee
where we have used  the divergence free condition and the
fact that
 \be
\int_{\mr^n}(h\cdot\nabla g)v\md x+\int_{\mr^n}(h\cdot\nabla
v)g\md x=0
 \ee
and
 \be
\int_{\mr^n}(w\cdot\nabla v)v\md x=0,\ \int_{\mr^n}(w\cdot\nabla
g)g\md x=0.
 \ee
Applying the trilinear estimates (\ref{2.22}) and (\ref{2.23}) in
Proposition \ref{prop2.1} to (\ref{3.64}),  it follows that
 \be
\|(v,g)\|^2_{L^2}+\int^t_0\|\nabla(v,g)\|^2_{L^2}\md s &\le&
\|(v_0,g_0)\|^2_{L^2}+\\ \nonumber &\ &
C\int^t_0\|(v,g)\|^2_{L^2}\|(u(s),b(s))\|^r_{\dot{B}^{\frac
np+\frac2r-1}_{p,r}}\md s.
 \ee
The Gronwall inequality yields the estimate (\ref{3.55}), and the
proof Theorem \ref{thm1.2} is thus complete.

\begin{rem}
The above arguments of formal computation can be justified by the
standard procedure of multiplying smoothing  sequence or frequency
localization operator (see \cite{G-P}).
\end{rem}


\section{Global well-posedness in Besov spaces: Case $n= 2$}
\setcounter{equation}{0}

 Now we are in position to prove of Theorem\ref{thm1.3}.
As in (\ref{1.13}) we decompose data $(u_0(x),b_0(x))$ into
$(v_0(x),g_0(x))\in L^2(\mr^2)$ and $(w_0(x),h_0(x))\in
\dot{B}^{2/\bar{p}-1}_{\bar{p},\bar{r}}(\mr^2)$ for some
$p<\bar{p}<\infty$ and $r<\bar{r}<\infty$ with small norm. By means of Theorem \ref{thm1.1},
we let $(w(t,x),h(t,x))\in C_b(\mr^+;\dot{B}^{2/{\bar p}-1}_{{\bar p},{\bar r}}(\mr^2))\cap
\widetilde{L}^{q}(\mr^+;\dot{B}^{\frac2p-1+2/q}_{{\bar p},{\bar r}}(\mr^2))$ solve globally the following system:
 \be\label{4.1add}
w_t-\tr w+\mb{P}\nabla\cdot(w\otimes w)-\mb{P}\nabla\cdot(h\otimes
h)=0,\\ \label{4.2add}
h_t-\tr h+\mb{P}\nabla\cdot(w\otimes
h)-\mb{P}\nabla\cdot(h\otimes
w)=0,\\
\label{4.3add}\mx{div} w=\mx{div} h=0,\\
\label{4.4add} w(x,0)=w_0(x),\ h(x,0)=h_0(x).
 \ee
and satisfies the regularity estimates  in Remark \ref{rem3.2}.
To prove Theorem\ref{thm1.3}, we are devoted to the study of the global well-posedness to a MHD-like system
 \be\nonumber
v_t-\tr v&+&(\mb{P}\nabla\cdot(v\otimes v)
+\mb{P}\nabla\cdot(v\otimes w)+\mb{P}\nabla\cdot(w\otimes v))\\
\label{4.1} &-&(\mb{P}\nabla\cdot(g\otimes
g)+\mb{P}\nabla\cdot(g\otimes h)+\mb{P}\nabla\cdot(h\otimes g))=0,
 \ee
 \be \nonumber
 g_t-\tr g&+&(\mb{P}\nabla\cdot(v\otimes g)+\mb{P}\nabla\cdot(v\otimes
h)+\mb{P}\nabla\cdot(w\otimes g))\\ \label{4.2}
&-&(\mb{P}\nabla\cdot(g\otimes v)+\mb{P}\nabla\cdot(h\otimes
v)+\mb{P}\nabla\cdot(g\otimes w))=0,
 \ee
 \be\label{4.3}
&\ & \mx{div}v=\mx{div}g=0,\\ \label{4.4} &\ & v(0,x)=v_0(x),\
g(0,x)=g_0(x).
 \ee
 for $L^2(\mr^2)$ data. To this goal, the first step is to establish the local well-posedness
of  the MHD-like equations (\ref{4.1})-(\ref{4.4})
for data $(v_0,g_0)\in
L^2(\mr^2)$ by means of  Gallagher-Planchon's argument\cite{G-P}.
 For completeness we give a clear presentation.  Let
 \be
X(I)={\cal C}_{\frac14}(I;L^4(\mr^2))\cap L^4(I;L^4(\mr^2))\cap
L^2(I;\dot{H}^1(\mr^2))\cap L^{2r,2}(I;L^{\frac{2r}{r-1}}(\mr^2))
 \ee
with norm
 \be\label{4.6}
\|f\|_X=\sup_{0<t<T}t^{1/4}\|f\|_{L^4}+\|f\|_{L^4(I;L^4)}+\|\nabla
f\|_{L^2(I;L^2)}+\|f\|_{L^{2r,2}(I;L^{\frac{2r}{r-1}})}.
 \ee
Here $I=[0,T]$ for some time $T>0$ and $r>1$. Then the local
well-posedness result is as follows.

\begin{thm}\label{thm4.1}
Let $(v_0,g_0)\in L^2(\mr^2)$. Then there exists a time $T>0$ and
the  unique solution $(v,g)\in X(I)$ to the system
(\ref{4.1})-(\ref{4.4}) such that
 \be
(v,g)\in C([0,T];L^2(\mr^2)).
 \ee
\end{thm}

\begin{proof}
For convenience we denote by $\|f\|_{X_1}$, $\|f\|_{X_2}$,
$\|f\|_{X_3}$ and $\|f\|_{X_4}$ every part of the norm
$\|\cdot\|_{X}$ in (\ref{4.6}). The MHD-like system
(\ref{4.1})-(\ref{4.4}) can be represented in the integral form
 \be\nonumber
v(t,x)&=&S(t)v_0-\int^t_0S(t-s)(\mb{P}\nabla\cdot(v\otimes v)
+\mb{P}\nabla\cdot(v\otimes w)+\mb{P}\nabla\cdot(w\otimes v))\md s\\
\label{4.7} &\ &+\int^t_0S(t-s)(\mb{P}\nabla\cdot(g\otimes
g)+\mb{P}\nabla\cdot(g\otimes h)+\mb{P}\nabla\cdot(h\otimes g))\md
s,
 \ee
 \be \nonumber
 g(t,x)&=&S(t)g_0-\int^t_0S(t-s)(\mb{P}\nabla\cdot(v\otimes g)+\mb{P}\nabla\cdot(v\otimes
h)+\mb{P}\nabla\cdot(w\otimes g))\md s\\ \label{4.8} &\ &
+\int^t_0S(t-s)(\mb{P}\nabla\cdot(g\otimes
v)+\mb{P}\nabla\cdot(h\otimes v)+\mb{P}\nabla\cdot(g\otimes w))\md
s.
 \ee
Here $S(t)$ denotes a semigroup whose kernel
$K_{\sqrt{t}}(x)=t^{-\frac n2}K(\frac x{\sqrt{t}})$, where
$K(\cdot)\in L^1(\mr^n)\cap L^\infty(\mr^n)$.
 We first consider the estimates of the free part $S(t)(v_0,g_0)$. Using Young's inequality and
Marcinkiewicz interpolation theorem \cite{M1} or \cite{S-W},  one
can easily prove that
 \be
\|S(t)(v_0,g_0)\|_{X_i}\le C\|(v_0,g_0)\|_{L^2}, \mx{ for } i=1,\
2,\ 3.
 \ee
For details please refer to \cite{Giga,M-G}. Choosing
$p=2r$ and $q=\frac{2r}{r-1}$ in Proposition \ref{prop2.2} we
obtain
 \be
\|S(t)(v_0,g_0)\|_{X_4}\le C\|(v_0,g_0)\|_{L^2}.
 \ee
It is necessary to point that $\|S(t)(v_0,g_0)\|_{X}$ can be small enough provided that
$T$ is small.

Next we estimate the bilinear terms. Observing that the bilinear terms
in integral equations (\ref{4.7})-(\ref{4.8}) are composed of two
kinds. One is the true nonlinear term $B(v,v)$, $B(g,g)$,
$B(g,v)$, $B(v,g)$; The other is $B(v,w)$, $B(w,v)$, $B(g,h)$,
$B(h,g)$, $B(v,h)$, $B(w,g)$, $B(h,v)$, $B(g,w)$, where $B(v,g)$
is a bilinear form with
 \be
B(v,g)=\int^t_0K_{\sqrt{t-s}}(x)*\mb{P}\nabla\cdot(v\otimes
g)(s)\md s,
 \ee
and $B(v,v)$, $B(g,v)$, $\cdots$ are similarly defined. We
consider the two kinds of bilinear forms, respectively.

A straightforward calculation yields
 \be\nonumber
\|B(v,g)\|_{L^4}(t)&\le&
C\int^t_0(t-s)^{-3/4}s^{-\frac14}\|g(s,\cdot)\|_{L^4}\md
s \|v\|_{X_1}\\ \label{4.12} &\le& C\int^t_0(t-s)^{-3/4}s^{-1/2}\md
s\|v\|_{X_1}\|g\|_{X_1}.
 \ee
Therefore one has
 \be
\|B(v,g)\|_{X_1}\le C \|v\|_{X_1}\|g\|_{X_1}.
 \ee
For the other kind of bilinear forms it can be shown  that
 \be\nonumber
\|B(v,w)\|_{L^4}(t)&\le&
C\int^t_0(t-s)^{-1/2}s^{-\frac12}\|v(s,\cdot)\|_{L^4}s^{\frac12}\|w(s,\cdot)\|_\infty \md s\\
&\le&\nonumber
C\ep_0\int^t_0(t-s)^{-1/2}s^{-\frac12}\|v(s,\cdot)\|_{L^4}\md s\\
\label{4.14} &\le& C\ep_0\int^t_0(t-s)^{-1/2}s^{-3/4}\md
s\|v\|_{X_1}.
 \ee
Here, we have used the estimate (\ref{globalE1}) for $w$
with $\alpha=0$ and $p=\infty$. Thus we obtain
 \be
\|B(v,w)\|_{X_1}\le C\ep_0\|v\|_{X_1}.
 \ee
 Applying the generalized
Young  inequality to the  second inequality in (\ref{4.12}) and the third inequality in (\ref{4.14})
with respect to time $t$,  respectively, we easily see that
 \be
\|B(v,g)\|_{X_2}\le C\|v\|_{X_1}\|g\|_{X_2},
 \ee
and
 \be
\|B(v,w)\|_{X_2}\le C\ep_0 \|v\|_{X_2}.
 \ee
Here, we have used the following  relations of indices
 \be
1+\frac14=\frac34+\frac14+\frac14,\quad
\frac14=\frac1\infty+\frac14+\frac1\infty,
 \ee
and
 \be
1+\frac14=\frac14+\frac12+\frac12, \quad
\frac14=\frac14+\frac1\infty+\frac1\infty.
 \ee
Similarly, to the norm $\|\cdot\|_{X_4}$ one has
 \be
\|B(v,g)\|_{X_4}\le C\|v\|_{X_4}\|g\|_{X_4}\mx{ or
}C\|v\|_{X_4}\|g\|_{X_1},
 \ee
and
 \be
\|B(v,w)\|_{X_4}\le C\ep_0\|v\|_{X_4}.
 \ee
Here, indices can be chosen as
 \be
1+\frac1{2r}=\frac1{2r}+\frac1{2r}+(1-\frac1{2r}), \quad
\frac12<\frac12+\frac12+\frac1\infty.
 \ee
or
 \be
1+\frac1{2r}=\frac34+\frac14+\frac1{2r}, \quad
\frac12=\frac1\infty+\frac1\infty+\frac12.
 \ee
and
 \be
1+\frac1{2r}=\frac12+\frac12+\frac1{2r}, \quad
\frac12=\frac1\infty+\frac1\infty+\frac12.
 \ee
Finally, we deal with the $X_3$ norm. A straightforward calculation shows
that
 \be\nonumber
\|B(\nabla v,g)\|_{L^2}(t)&\le& C\int^t_0(t-s)^{-3/4}\|\nabla
v(s,\cdot)\|_{L^2}\|g(s,\cdot)\|_{L^4}\md s\\ &\le&
C\int^t_0(t-s)^{-3/4}s^{-\frac14}\|\nabla v(s,\cdot)\|_{L^2}\md s\|g\|_{X_1}.
 \label{4.32}\ee
Applying the generalized Young  inequality to (\ref{4.32}),
 one has
 \be
\|B(v,g)\|_{X_3}\le
C(\|v\|_{X_3}\|g\|_{X_1}+\|g\|_{X_3}\|v\|_{X_1}).
 \ee
Similar computation follows that
 \be
\|B(\nabla v,w)\|_{L^2}\le C\int^t_0(t-s)^{-\frac1r-\frac12}s^{-(\frac12-\frac1r)}\|\nabla
v(s,\cdot)\|_{L^2}\md s \|w(s)\|_{{\cal C}_{\frac12-\frac1r}(I; L^r)},
 \ee
and
 \be
\|B(v,\nabla w)\|_{L^2}\le
C\int^t_0(t-s)^{-\frac1{2r}-\frac12}s^{-(1-\frac1r)}\|v(s,\cdot)\|_{L^{\frac{2r}{r-1}}}\md s
\|\nabla w(s)\|_{{\cal C}_{1-\frac1r}(I; L^r)}\md s.
 \ee
Using the generalized Young inequality and
applying the estimate (\ref{globalE1}) to $w$ with $\alpha=0$ and
$\alpha=1$, respectively, we arrive at
 \be
\|B(v,w)\|_{X_3}\le C\ep_0(\|v\|_{X_3}+\|v\|_{X_4}).
 \ee
Here, the indices chosen satisfy
 \be
1+\frac12=(\frac1r+\frac12)+(\frac12-\frac1r)+\frac12, \quad
\frac12=\frac1\infty+\frac1\infty+\frac12,
 \ee
and
 \be
1+\frac12=(\frac1{2r}+\frac12)+(1-\frac1r)+\frac1{2r}, \quad
\frac12=\frac1\infty+\frac1\infty+\frac12,
 \ee
respectively.

 In conclusion, we obtain the two kinds of bilinear estimates as
 \be\label{4.32}
\|B(v,g)\|_X\le C\|v\|_X\|g\|_X,\  \|B(v,w)\|_X\le C\ep_0\|v\|_X,
 \ee
for some constant C.

Assume $(v_1,g_1)$ and $(v_2,g_2)$ are two solutions (\ref{4.7})
and (\ref{4.8}), then we have the difference
 \be
v_1-v_2&=&B(g_1-g_2,g_1)+G(g_2,g_1-g_2)+B(g_1-g_2,h)+G(h,g_1-g_2)\\\nonumber
&\ & - (B(v_1-v_2,v_1)+G(v_2,v_1-v_2)+B(v_1-v_2,w)+G(w,v_1-v_2)),
 \ee
and
 \be
g_1-g_2&=&B(g_1-g_2,v_1)+G(g_2,v_1-v_2)+B(h,v_1-v_2)+G(g_1-g_2,w)\\\nonumber
&\ & - (B(v_1-v_2,g_1)+G(v_2,g_1-g_2)+B(v_1-v_2,h)+G(w,g_1-g_2)).
 \ee
Arguing similarly as in deriving (\ref{4.32}) one has
 \be
\|(v_1-v_2,g_1-g_2)\|_X\le
C\|(v_1-v_2,g_1-g_2)\|_X(\|(v_1,g_1)\|_X+\|(v_2,g_2)\|_X+\ep_0).
 \ee
We can construct a complete metric space as following:
 \be
E=\Big\{(f, g):\;  (f, g)\in X(I);\ \|(f, g)\|_X\le 2K_0(T)\Big\}
 \ee
with metric $d((f_1,f_2), (g_1,g_2))=\|(f_1-f_2, g_1-g_2)\|_X$, where
$K_0(T)=\|S(t)(v_0,g_0)\|_X$ can be small enough provided we
choose the existent time $T$ small. Thus Applying Banach contraction
mapping principle to nonlinear mapping ${\cal T}$ defined by
the right sides of (\ref{4.7})-(\ref{4.8}) on $E$, a unique local solution
$(v,g)$ to (\ref{4.7})-(\ref{4.8}) is obtained on interval $[0,T]$.
 Using the fact $(v,g)\in X(I)$
 it is not difficult to verify
that $(v,g)\in L^\infty([0,T];L^2(\mr^2))\cap
L^2([0,T];\dot{H}^1(\mr^2))$. Noting the imbedding relation
$L^2(\mr^2)\hookrightarrow \dot{B}^{2/p-1}_{p,q}(\mr^2)$ for $p>2$
and $q\ge 2$,  one also has $(v,g)\in
C([0,T];\dot{B}^{2/p-1}_{p,q}(\mr^2))$. The proof of Theorem
\ref{thm4.1} is thus complete.
\end{proof}

To complete the proof of global existence of solution we need the
following energy inequality.

\begin{lem}\label{lem4.1}
There exists a time $t_0$ and constant $C(\ep_0,t_0)$ such that
the solution $(v,g)$ constructed in Theorem \ref{thm4.1} satisfies
the inequality
 \be
\sup_{t_0<s<t}\|(v(s),g(s))\|^2_{L^2}+\|\nabla(v,g)\|^2_{L^2([t_0,t];L^2)}\le
2C\|(v_0,g_0)\|^2_{L^2}\Big(\frac t{t_0}\Big)^{\ep_0}.
 \ee
 In particular, if $\frac 2p+\frac2r>1$, one has
 \be\label{4.39}
\|(v(s),g(s))\|^2_{L^2}+\int^t_0\|\nabla(v,g)\|^2_{L^2}\md s\le
C(\ep_0, t_0)\|(v_0,g_0)\|^2_{L^2}.
 \ee
\end{lem}

\begin{proof}
Formally, we multiply (\ref{4.1}) by $v(t,x)$, (\ref{4.2}) by
$g(t,x)$, and integrate both equations with respect to $x$ on
$\mr^2$ and $t$ on $[0,T]$, then sum the results to arrive at
 \be\nonumber
&\ &
\|(v(t),g(t))\|^2_{L^2}+2\int^t_0\|\nabla(v(s),g(s))\|^2_{L^2}\md
s \le 2\|(v_0,g_0)\|^2_{L^2}\\ \nonumber &\ &
-2\int^t_0\int_{\mr^2}(v\cdot\nabla v)w(s, x)\md x\md
s-2\int^t_0\int_{\mr^2}(g\cdot\nabla g)w(s, x)\md x\md
s\\\label{4.38} &\ & -2\int^t_0\int_{\mr^2}(g\cdot\nabla
v)h(s, x)\md x\md s+2\int^t_0\int_{\mr^2}(v\cdot\nabla g)h(s, x)\md
x\md s,
 \ee
where we have made use of  the fact
 \be
\int_{\mr^2}(g\cdot \nabla g)v(t,x)\md x+\int_{\mr^2}(g\cdot
\nabla v)g(t,x)\md x=0,\\ \int_{\mr^2}(h\cdot \nabla g)v(t,x)\md
x+\int_{\mr^2}(h\cdot \nabla v)g(t,x)\md x=0,
 \ee
and the divergence free condition. For simplicity of notation,
$\|(v(t),g(t))\|^2_{L^2}$ denotes
$\|v(t)\|^2_{L^2}+\|g(t)\|^2_{L^2}$ and
$\|(w,h)\|^r_{\dot{B}^{2/p+2/r-1}_{p,r}}$ denotes
$\|w\|^r_{\dot{B}^{2/p+2/r-1}_{p,r}}+\|h\|^r_{\dot{B}^{2/p+2/r-1}_{p,r}}$.
Applying estimates (\ref{2.22})-(\ref{2.23}) in Proposition
\ref{prop2.1} to (\ref{4.38}) it follows that
 \be
\|(v(t),g(t))\|^2_{L^2}+\int^t_0\|\nabla(v(s),g(s))\|^2_{L^2}\md s
&\le& 2\|(v_0,g_0)\|^2_{L^2}\\ \nonumber &\ &
+C\int^t_0\|(v,g)\|^2_{L^2}\|(w,h)\|^r_{\dot{B}^{2/p+2/r-1}_{p,r}}\md
s.
 \ee
Noting that
$\int^\infty_0\|(w,h)\|^r_{\dot{B}^{2/p+2/r-1}_{p,r}}\md s\le
\ep_0$ by Theorem \ref{thm1.1} and applying the Gronwall
inequality yields the estimate (\ref{4.39}). The above formal
arguments can be again justified by standard procedure, see
\cite{G-P}.

 According to Theorem \ref{thm4.1} there exists a local solution
on $[0,T)$ and there exists a small time $t_0$ such that the local
solution is smoothed out on $[t_0,T)$. So the energy estimate
(\ref{4.38}) is justified if we replace $0$ with $t_0$. Using the
fact $\sup_{0<t<\infty}\sqrt{t}\|(w(t),h(t)\|_{L^\infty}<\ep_0$
and H\"older and Young inequalities one has
 \be
\int^t_{t_0}\int_{\mr^2}(v\cdot\nabla v)w(s, x)\md x\md s\le
\ep_0\bigg(\int^t_{t_0}\|\nabla v(s)\|^2_{L^2}\md
s+\frac14\int^t_{t_0}\frac{\|v(s)\|^2_{L^2}}s\md s\bigg).
 \ee
Similarly, we can estimate $\int^t_{t_0}\int_{\mr^2}(g\cdot\nabla
g)w(s, x)\md x\md s$, $\int^t_0\int_{\mr^2}(g\cdot\nabla
v)h(s, x)\md x\md s$ and $\int^t_{t_0}\int_{\mr^2}(v\cdot\nabla
g)h(s, x)\md x\md s$. Inserting them to the estimate (\ref{4.38})
we  arrive at
 \be\nonumber
&\ &
\|(v(t),g(t))\|^2_{L^2}+2\int^t_0\|\nabla(v(s),g(s))\|^2_{L^2}\md
s \le 2\|(v_0,g_0)\|^2_{L^2}\\ &\ &
+\ep_0\bigg(4\int^t_{t_0}\|\nabla (v(s),g(s))\|^2_{L^2}\md
s+\int^t_{t_0}\frac{\|(v(s),g(s))\|^2_{L^2}}s\md s\bigg).
 \ee
Applying the Gronwall inequality we thus complete the proof of Lemma
\ref{lem4.1}.
\end{proof}
By Lemma \ref{lem4.1} we obtain that there exists a global
solution $(v(t,x),g(t,x))$ to (\ref{4.1})-(\ref{4.4}) for the
$L^2$ data $(v_0(x),g_0(x))$. This together with Theorem \ref{thm1.1}
yields Theorem \ref{thm1.3} except for the
estimate (\ref{1.9}).

To prove the estimate (\ref{1.9}) we first recall  the concept
of the real interpolation method \cite{B-L}. Let $X_1$, $X_2$ and
$X$ be three Banach spaces, $X$ is the real interpolation space of
$X_1$ and $X_2$ with
 \be
X=[X_1,X_2]_{\theta,r} \mx{ with } 0\le\theta\le 1 \mx{ and } 1\le
r\le \infty.
 \ee
Then for any $f\in X$ one has
 \be
\|f\|_X=\bigg(\sum_{j\in \mb{Z}}2^{jr\theta}K(f,j)^r\bigg)^{1/r},
 \ee
where $K(f,j)=\inf_{g\in X_2}(\|f-g\|_{X_1}+2^{-j}\|g\|_{X_2})$.
In our case of Theorem \ref{thm1.3} we take
$X=\dot{B}^{2/p-1}_{p,r}(\mr^2)$,
$X_1=\dot{B}^{2/\bar{p}-1}_{\bar{p},\bar{r}}(\mr^2)$ and
$X_2=L^2(\mr^2)$, where $p<\bar{p}$ and $r<\bar{r}$. The following
lemma states the property of the interpolation norm, see
\cite{G-P} for details.

\begin{lem}\label{lem4.2}
Let $X$ be the real interpolation space of two Banach spaces $X_1$
and $X_2$ with relation $X_2\hookrightarrow X\hookrightarrow X_1$,
then there exists a constant $C$ such that for any integer $j_0\ge
1 $ and $f\in X$, one has
 \be\label{4.48}
\bigg(\sum_{j\ge j_0}2^{jr\theta}K(f,j)^r\bigg)^{1/r}\le
\|f\|_X\le C2^{j_0}\bigg(\sum_{j\ge
j_0}2^{jr\theta}K(f,j)^r\bigg)^{1/r}.
 \ee
\end{lem}

For any $(u_0,b_0)\in X$,  one may decompose $(u_0,b_0)$ as
 \be
(u_0,b_0)=(v_0,g_0)+(w_0,h_0),
 \ee
with $(v_0,g_0)\in X_2$, $(w_0,h_0)\in X_1$ and
$\|(w_0,h_0)\|_{X_1}\le\ep_0$. The associated solution is
$(u,b)=(v,g)+(w,h)$, where $(v,g)$ and $(w,h)$ are the solutions
constructed in Theorem \ref{thm4.1} and \ref{thm1.1} with initial data $(w_0,h_0)$,
respectively. Noting (\ref{4.39}) in Lemma \ref{lem4.1} one has
the a priori estimates
 \be \label{4.50}
\|(v,g)\|_{X_2}\le C(\ep_0)\|(v_0,g_0)\|_{X_2}, \mx{ and }
\|(w,h)\|_{X_1}\le 2\|(w_0,h_0)\|_{X_1}.
 \ee

The first inequality  of (\ref{4.48}) in Lemma \ref{lem4.2} implies
that there exist $(w_0^j,h_0^j)\in X_1$ and $(v_0^j,g_0^j)\in X_2$
such that
 \be\label{4.51}
\|(u_0,b_0)\|_{X}\ge
\frac12\|2^{j\theta}(\|(w_0^j,h_0^j)\|_{X_1}+2^{-j}\|(v_0^j,g_0^j)\|_{X_2})\|_{l^r(j\ge
j_0)}.
 \ee
Thus for any $j\ge j_0$ one has
 \be
\|(w_0^j,h_0^j)\|_{X_1}\le 2\cdot2^{-j_0\theta}\|(u_0,b_0)\|_X<\ep_0,
 \ee
where we  choose $j_0$ such that $2^{-j_0\theta }=\frac
{\ep_0}{2\|(u_0,b_0)\|_X}$.

We construct solutions $(v^j,g^j)$ and $(w^j,h^j)$ to the Cauchy
problem (\ref{4.1})-(\ref{4.4}) and (\ref{3.1})-(\ref{3.4})
associated with data $(v_0^j,g_0^j)$ and $(w_0^j,h_0^j)$, together
with (\ref{4.50}) and (\ref{4.51}).  It follows that
 \be
\|(u_0,b_0)\|_X\ge
C^{-1}\|2^{j\theta}(\|(w^j,h^j)\|_{X_1}+2^{-j}\|(v^j,g^j)\|_{X_2})\|_{l^r(j\ge
j_0 )}.
 \ee
In view of the definition of functional $K(f,j)$ and
$(u^j,b^j)=(v^j,g^j)+(w^j,h^j)$ we derive that
 \be
\|(u_0,b_0)\|_X\ge C^{-1}\|2^{j\theta}K((u,b),j)\|_{l^r(j\ge j_0)}.
 \ee
Noting the right estimate (\ref{4.48}) in Lemma \ref{lem4.2} and
the choice of $j_0$ we finally arrive at
 \be
\|(u,b)\|_{X}\le C(\theta,r)\|(u_0,b_0)\|_X^{1+1/\theta}.
 \ee
Here $\frac1\theta > \frac p2$, $\theta=\frac
{2(\bar{p}-p)}{p(\bar{p}-2)}$ is the interpolation parameter and
$C(\theta,r)\le C\frac {2^\theta}{(2^{\theta r}-1)^{1/r}}\le
C(p,r)$. Choosing $\beta=\frac 1\theta$ the estimate (\ref{1.9})
follows.

\section*{Acknowledgements}  The authors thank the referees and the
associated editor for their invaluable comments and suggestions
which helped improve the paper greatly. C. Miao  was partly
supported by the NSF of China (No.10725102). B.Yuan was partially
supported by the NSF of China (No. 10771052), the NSF of Henan
Province (No. 0611055500) and the Doctoral Foundation of Henan
Polytechnic University.


\end{document}